\documentclass[12pt,a4paper]{amsart}

\usepackage[T1]{fontenc}
\usepackage[utf8]{inputenc}
\usepackage{lmodern}
\usepackage{amsmath,amssymb,mathtools}
\usepackage{geometry}
\usepackage{microtype}
\usepackage{enumitem}
\usepackage{xcolor}
\usepackage[colorlinks=true,linkcolor=blue!55!black,citecolor=blue!55!black,urlcolor=blue!55!black]{hyperref}

\setlist{nosep}
\numberwithin{equation}{section}

\newtheorem{theorem}{Theorem}[section]
\newtheorem{proposition}[theorem]{Proposition}
\newtheorem{lemma}[theorem]{Lemma}
\newtheorem{corollary}[theorem]{Corollary}
\theoremstyle{remark}

\newcommand{\D}{\mathbb D}

\title{Pentablock and Hexablock are Lempert Domains}
\author{Armen Edigarian}


\begin{document}

\begin{abstract}
We prove that the hexablock admits an explicit compactly nested exhaustion by bounded domains with real-analytic strongly linearly convex boundary.   The lower-dimensional classical domains: the tetrablock, pentablock, and symmetrized bidisc occur as complex linear sections of the hexablock, and the same level sets restrict to strongly linearly convex exhaustions of all three.  Consequently the hexablock and the pentablock are Lempert domains; the same argument recovers the known Lempert property of the tetrablock and the symmetrized bidisc.  We also give a direct continuous negative plurisubharmonic exhaustion of the hexablock, hence all four domains are hyperconvex and taut.
\end{abstract}
\maketitle

\section{Introduction}
\label{sec:introduction}

\subsection{Invariant distances and Lempert theory}
Let $\rho$ denote the Poincar\'e distance on the unit disc $\D$ and let $D\subset\mathbb C^n$ be a domain.  The
Carath\'eodory distance and the Lempert function are
\begin{align}
 c_D(z,w)
 &=\sup\{\rho(f(z),f(w)):f\in\mathcal O(D,\mathbb D)\},\label{eq:cD}\\
 \ell_D(z,w)
 &=\inf\{\rho(\alpha,\beta):f\in\mathcal O(\mathbb D,D),
                    f(\alpha)=z,\ f(\beta)=w\}.\label{eq:lD}
\end{align}
The Kobayashi distance $k_D$ is the chain envelope of $\ell_D$, and hence
\[
 c_D\le k_D\le \ell_D.
\]
We denote the infinitesimal Carath\'eodory--Reiffen and Kobayashi--Royden
metrics by $\gamma_D$ and $\kappa_D$, respectively.  We refer to
Jarnicki--Pflug \cite{JP2013} for general background.

Lempert's theorem \cite{Lempert1981,Lempert1984} gives
\[
 c_D=k_D=\ell_D,\qquad \gamma_D=\kappa_D
\]
for bounded convex domains.  For $n\ge2$, the same conclusion holds on
bounded strongly linearly convex domains in $\mathbb C^n$ with $C^2$ boundary;
in this form we use Kosi\'nski--Warszawski \cite{KW2013}.  We shall call a
taut domain a \emph{Lempert domain} when $c_D=\ell_D$.

\subsection{Four domains arising from \texorpdfstring{$\mu$}{mu}-synthesis}

For $\lambda=(\lambda_1,\lambda_2)\in\mathbb D^2$, let
\[
 \pi(\lambda)=(\lambda_1+\lambda_2,\lambda_1\lambda_2).
\]
The symmetrized bidisc is
\[
 \mathbb G_2=\pi(\mathbb D^2).
\]
It arose in the spectral Nevanlinna--Pick problem studied by Agler and Young
\cite{AY1999,AY2000}.  Agler and Young established the Lempert property of
$\mathbb G_2$ in their analysis of its hyperbolic geometry \cite{AY2004}.
Costara used their work to show, in particular, that $\mathbb G_2$ is not
biholomorphic to a convex domain \cite{Costara2004}; Pflug and Zwonek later
gave a geometric proof of the Lempert equality by exhausting $\mathbb G_2$
by strongly linearly convex domains \cite{PZ2012}.

We write points of the tetrablock consistently as
\[
 x=(x_1,x_2,x_3)\in\mathbb C^3.
\]
The tetrablock is
\[
 \mathbb E
 =\left\{x:\ 1-x_1\lambda_1-x_2\lambda_2
       +x_3\lambda_1\lambda_2\ne0
       \text{ for every }\lambda=(\lambda_1,\lambda_2)
       \in\overline{\mathbb D}^{\,2}\right\}.
\]
It was introduced by Abouhajar, White and Young \cite{AWY2007}.  The
Lempert theorem for $\mathbb E$ was proved in
Edigarian--Kosi\'nski--Zwonek \cite{EKZ2013}; its $\mathbb C$-convexity was
proved by Zwonek \cite{Zwonek2013}.

The pentablock, introduced by Agler, Lykova and Young \cite{ALY2015}, is
\[
 \mathbb P
 =\{(a,s,p)=(A_{21},\operatorname{tr}A,\det A):
                    A\in M_2(\mathbb C),\ \|A\|<1\}.
\]
It is a Hartogs domain over $\mathbb G_2$.  Its automorphism group was
determined by Kosi\'nski \cite{Kosinski2015}; Zapa\l{}owski proved that it
cannot be exhausted by domains biholomorphic to convex domains
\cite{Zapalowski2015}; and Su proved that it is $\mathbb C$-convex
\cite{Su2020}.

The hexablock $\mathbb H\subset\mathbb C^4$ was introduced by
Biswas--Pal--Tomar \cite{BPT2025}.  It is a Hartogs domain over the
tetrablock and contains canonical copies of $\mathbb E$, $\mathbb P$, and
$\mathbb G_2$; in fact these domains occur as analytic retracts of
$\mathbb H$ \cite{BPT2026}.  The point of the present paper is that the
strongly linearly convex approximation can be proved once, on $\mathbb H$,
and then inherited by the three classical sections.

\subsection{Main result}

For $x=(x_1,x_2,x_3)\in\mathbb E$ set
\[
 q(x)=x_1x_2-x_3,\qquad
 K(x)=1-|x_1|^2-|x_2|^2+|x_3|^2,
\]
\[
 \Delta(x)=K(x)^2-4|q(x)|^2,
 \qquad
 V(x)=\frac{K(x)+\sqrt{\Delta(x)}}2,
\]
and
\[
 F(a,x)=\sqrt{\Delta(x)}\bigl(V(x)-|a|^2\bigr).
\]
The technical result of the paper is the following.

\begin{theorem}\label{thm:intro-main}
For every $0<\varepsilon<1$,
\[
 D_\varepsilon=\{(a,x)\in\mathbb H:F(a,x)>\varepsilon\}
\]
is a bounded connected domain with real-analytic strongly linearly convex
boundary.  Moreover,
\[
 \overline{D_\varepsilon}\subset D_{\varepsilon'}\Subset\mathbb H
 \quad(0<\varepsilon'<\varepsilon<1),
 \qquad
 \bigcup_{\varepsilon>0}D_\varepsilon=\mathbb H.
\]
Thus $\mathbb H$ admits a compactly nested strongly linearly convex
exhaustion.
\end{theorem}

The three relevant complex linear sections are
\[
 S_{\mathbb E}=\{a=0\},\qquad
 S_{\mathbb P}=\{x_1=x_2\},\qquad
 S_{\mathbb G}=\{a=0,\ x_1=x_2\}.
\]
They satisfy
\[
 \mathbb H\cap S_{\mathbb E}\simeq\mathbb E,
 \qquad
 \mathbb H\cap S_{\mathbb P}\simeq\mathbb P,
 \qquad
 \mathbb H\cap S_{\mathbb G}\simeq\mathbb G_2.
\]
The restrictions $D_\varepsilon\cap S$ are again regular, connected and
strongly linearly convex.  No second convexity calculation is required.
Consequently we obtain the following statement from the single exhaustion
above.

\begin{corollary}\label{cor:intro-domains}
Each of
\[
 \mathbb H,\qquad \mathbb P,\qquad \mathbb E,\qquad \mathbb G_2
\]
admits an increasing compactly nested exhaustion by bounded domains with
real-analytic strongly linearly convex boundary.  All four domains are
hyperconvex and taut, and
\[
 c_D=k_D=\ell_D,\qquad \gamma_D=\kappa_D
\]
for $D=\mathbb H,\mathbb P,\mathbb E,\mathbb G_2$.
In particular, the hexablock and the pentablock are Lempert domains.
\end{corollary}

For $\mathbb G_2$ this recovers the exhaustion theorem of Pflug--Zwonek,
and for $\mathbb E$ it gives another proof of the already known Lempert
equality.  The novelty is that the same four-dimensional level function
also settles the hexablock and the pentablock.

The paper is organized as follows.  Sections~\ref{sec:hexablock-exhaustion}--
\ref{sec:regularity-exhaustion} prove Theorem~\ref{thm:intro-main}.  Section
\ref{sec:hyperconvexity} proves tautness by a direct plurisubharmonic
argument.  Section~\ref{sec:sections} identifies the three linear sections
and transfers the exhaustion to them.  The final section passes to the
invariant distances and metrics.

\section{The hexablock and the explicit exhaustion}
\label{sec:hexablock-exhaustion}

Write $x=(x_1,x_2,x_3)$ and $\lambda=(\lambda_1,\lambda_2)$.  Put
\[
 h_x(\lambda_1,\lambda_2)
 =1-x_1\lambda_1-x_2\lambda_2+x_3\lambda_1\lambda_2.
\]
and let
\[
 \mathbb E=\{x\in\mathbb C^3:h_x(\lambda_1,\lambda_2)\ne0
       \text{ for every }\lambda=(\lambda_1,\lambda_2)\in\overline{\mathbb D}^{\,2}\}
\]
be the tetrablock. We use the hexablock of Biswas--Pal--Tomar
\cite{BPT2025}, namely
\begin{equation}\label{hex:def}
 \mathbb H=\left\{(a,x)\in\mathbb C\times\mathbb E:
 |a|\sup_{\lambda=(\lambda_1,\lambda_2)\in\mathbb D^2}
 \frac{\sqrt{(1-|\lambda_1|^2)(1-|\lambda_2|^2)}}{|h_x(\lambda_1,\lambda_2)|}<1\right\}.
\end{equation}
Throughout, $\mathbb H$ denotes this Hartogs domain.

For $x=(x_1,x_2,x_3)\in\mathbb E$ set
\begin{align}
 q(x)&=x_1x_2-x_3,\label{hex:q}\\
 K(x)&=1-|x_1|^2-|x_2|^2+|x_3|^2,\label{hex:K}\\
 \Delta(x)&=K(x)^2-4|q(x)|^2,\label{hex:Delta}\\
 V(x)&=\frac{K(x)+\sqrt{\Delta(x)}}2.\label{hex:V}
\end{align}
The square root in \eqref{hex:V} will always be the positive square root.
The proof below establishes that $K>0$ and $\Delta>0$ on $\mathbb E$.
Define
\begin{equation}\label{hex:F}
 F(a,x)=\sqrt{\Delta(x)}\bigl(V(x)-|a|^2\bigr).
\end{equation}

\begin{theorem}\label{hex:exhaustion-theorem}
The hexablock has the explicit description
\[
 \mathbb H=\{(a,x):x\in\mathbb E,\ |a|^2<V(x)\}.
\]
For every $0<\varepsilon<1$, the set
\begin{equation}\label{hex:D-epsilon}
 D_\varepsilon=\{(a,x)\in\mathbb H:F(a,x)>\varepsilon\}
\end{equation}
is a bounded connected domain with real-analytic strongly linearly convex
boundary. Equivalently, with $K,q,\Delta$ as above,
\begin{equation}\label{hex:fully-explicit-exhaustion}
 \begin{split}
 D_\varepsilon=\{(a,x_1,x_2,x_3):\;&|x_3|<1,\quad K>2|q|,\\
 &\Delta+(K-2|a|^2)\sqrt\Delta>2\varepsilon\}.
 \end{split}
\end{equation}
If $0<\varepsilon'<\varepsilon<1$, then
\[
 \overline{D_\varepsilon}\subset D_{\varepsilon'}\Subset\mathbb H,
 \qquad
 \bigcup_{0<\varepsilon<1}D_\varepsilon=\mathbb H.
\]
Thus, for example, $D_{2^{-n}}$, $n\ge1$, is a nested strongly linearly
convex exhaustion of $\mathbb H$.
\end{theorem}

The defining function in \eqref{hex:D-epsilon} combines both possible
boundary degenerations. The factor $V-|a|^2$ vanishes at the boundary of a
Hartogs fibre, whereas $\sqrt\Delta$ vanishes when the tetrablock base
approaches its boundary. The proof of strong linear convexity uses the
projective covariance of their product.

\section{Projective automorphisms and normal form}

For $\alpha,\beta\in\mathbb D$ put
\[
 c=(1-|\alpha|^2)(1-|\beta|^2),
 \qquad
 d_{\alpha,\beta}(x)
 =1-\bar\alpha x_1-\bar\beta x_2+\bar\alpha\bar\beta x_3.
\]
The denominator does not vanish on $\mathbb E$. Define
\begin{equation}\label{hex:projective-base-map}
 \begin{split}
 \Phi_{\alpha,\beta}(x)
 =\frac1{d_{\alpha,\beta}(x)}\bigl(
 &x_1-\alpha-\bar\beta x_3+\alpha\bar\beta x_2,\\
 &x_2-\beta-\bar\alpha x_3+\beta\bar\alpha x_1,\\
 &x_3-\beta x_1-\alpha x_2+\alpha\beta\bigr),
 \end{split}
\end{equation}
and
\begin{equation}\label{hex:projective-full-map}
 \Psi_{\alpha,\beta}(a,x)
 =\left(\frac{\sqrt c\,a}{d_{\alpha,\beta}(x)},
           \Phi_{\alpha,\beta}(x)\right).
\end{equation}
These are projective linear maps in the affine coordinates under
consideration.

\begin{lemma}\label{hex:automorphisms-normal-form}
The maps $\Phi_{\alpha,\beta}$ and $\Psi_{\alpha,\beta}$ are automorphisms
of $\mathbb E$ and $\mathbb H$, respectively, with inverses obtained by
replacing $(\alpha,\beta)$ by $(-\alpha,-\beta)$. Every point of
$\mathbb H$ can be carried, by one such map followed by coordinate
rotations, to
\[
 (A,0,0,r),\qquad 0\le A<1,\quad 0\le r<1.
\]
Moreover,
\begin{align}
 q(\Phi_{\alpha,\beta}(x))
 &=\frac{c\,q(x)}{d_{\alpha,\beta}(x)^2},\label{hex:q-covariance}\\
 K(\Phi_{\alpha,\beta}(x))
 &=\frac{c\,K(x)}{|d_{\alpha,\beta}(x)|^2},\label{hex:K-covariance}\\
 \Delta(\Phi_{\alpha,\beta}(x))
 &=\frac{c^2\Delta(x)}{|d_{\alpha,\beta}(x)|^4}.
 \label{hex:Delta-covariance}
\end{align}
\end{lemma}

\begin{proof}
For $\gamma\in\mathbb D$ define the disc automorphism
\[
 \tau_\gamma(\zeta)=\frac{\bar\gamma+\zeta}{1+\gamma\zeta},\qquad \zeta\in\mathbb D.
\]
A direct expansion yields
\begin{equation}\label{hex:kernel-transform}
 d_{\alpha,\beta}(x)h_{\Phi_{\alpha,\beta}(x)}(\lambda_1,\lambda_2)
 =(1+\alpha \lambda_1)(1+\beta \lambda_2)
 h_x\bigl(\tau_\alpha(\lambda_1),\tau_\beta(\lambda_2)\bigr).
\end{equation}
Since $\tau_\alpha$ and $\tau_\beta$ are automorphisms of the disc, this
identity proves that $\Phi_{\alpha,\beta}$ maps $\mathbb E$ into itself.
Its inverse formula follows either by substitution or by applying the
inverse disc automorphisms in \eqref{hex:kernel-transform}.

Also,
\[
 1-|\tau_\alpha(\lambda_1)|^2
 =\frac{(1-|\alpha|^2)(1-|\lambda_1|^2)}{|1+\alpha \lambda_1|^2}.
\]
Together with \eqref{hex:kernel-transform}, this gives, if
$(a',x')=\Psi_{\alpha,\beta}(a,x)$,
\[
 \frac{|a'|\sqrt{(1-|\lambda_1|^2)(1-|\lambda_2|^2)}}{|h_{x'}(\lambda_1,\lambda_2)|}
 =\frac{|a|\sqrt{(1-|\tau_\alpha(\lambda_1)|^2)
                  (1-|\tau_\beta(\lambda_2)|^2)}}
 {|h_x(\tau_\alpha(\lambda_1),\tau_\beta(\lambda_2))|}.
\]
Taking suprema proves the assertion about $\mathbb H$ and its inverse.

For completeness, the covariance identities have a short matrix proof.
Set
\[
 M(x)=\begin{pmatrix}1&-x_2\\-x_1&x_3\end{pmatrix},
 \qquad J=\begin{pmatrix}1&0\\0&-1\end{pmatrix},
\]
so that
\[
 \det M(x)=-q(x),\qquad
 K(x)=\operatorname{tr}\bigl(JM(x)JM(x)^*\bigr).
\]
For
\[
 L_\alpha=\begin{pmatrix}1&\bar\alpha\\\alpha&1\end{pmatrix},
 \qquad
 R_\beta=\begin{pmatrix}1&\beta\\\bar\beta&1\end{pmatrix}
\]
one has
\[
 M(\Phi_{\alpha,\beta}(x))
 =\frac{L_\alpha M(x)R_\beta}{d_{\alpha,\beta}(x)},
\]
\[
 L_\alpha^*JL_\alpha=(1-|\alpha|^2)J,
 \qquad R_\beta JR_\beta^*=(1-|\beta|^2)J.
\]
Taking determinants and using cyclicity of the trace proves
\eqref{hex:q-covariance} and \eqref{hex:K-covariance};
\eqref{hex:Delta-covariance} follows.

It remains to prove the normal form. For fixed $x\in\mathbb E$, the
function
\[
 \kappa_x(\lambda_1,\lambda_2)=
 \frac{\sqrt{(1-|\lambda_1|^2)(1-|\lambda_2|^2)}}{|h_x(\lambda_1,\lambda_2)|}
\]
is continuous on $\overline{\mathbb D}^2$, vanishes on its boundary,
and is positive at $(0,0)$. It therefore attains its maximum at some
$(\lambda_{1,0},\lambda_{2,0})\in\mathbb D^2$. Differentiating $\log\kappa_x$ at this point
and using $h_x\ne0$ gives
\[
 \bar \lambda_{1,0}=\frac{x_1-x_3\lambda_{2,0}}{1-x_2\lambda_{2,0}},
 \qquad
 \bar \lambda_{2,0}=\frac{x_2-x_3\lambda_{1,0}}{1-x_1\lambda_{1,0}}.
\]
Choose $\alpha=\bar \lambda_{1,0}$ and $\beta=\bar \lambda_{2,0}$. These identities say
exactly that the first two numerators in
\eqref{hex:projective-base-map} vanish. Hence
\[
 \Phi_{\alpha,\beta}(x)=(0,0,r_0),\qquad |r_0|<1.
\]
The last inequality follows directly from membership in $\mathbb E$.
The rotations
\[
 (a,x_1,x_2,x_3)\longmapsto
 (\eta a,\omega_1x_1,\omega_2x_2,\omega_1\omega_2x_3),
 \qquad \eta,\omega_1,\omega_2\in\partial\mathbb D,
\]
preserve \eqref{hex:def}; they make the first and last coordinates
nonnegative real. At $(0,0,r_0)$ the supremum of $\kappa$ is $1$, as
proved below, so the first coordinate has modulus less than $1$.
\end{proof}

\section{The radius and the boundary behaviour of \texorpdfstring{$F$}{F}}

\begin{lemma}\label{hex:radius-and-boundary}
One has $K>0$ and $\Delta>0$ on $\mathbb E$, and the radius of the
hexablock fibre over $x$ is $\sqrt{V(x)}$. The function $F$ is positive
and real-analytic on $\mathbb H$, extends continuously to
$\overline{\mathbb H}$, and vanishes on $\partial\mathbb H$. It satisfies
\begin{equation}\label{hex:F-covariance}
 F(\Psi_{\alpha,\beta}(a,x))
 =\frac{c^2F(a,x)}{|d_{\alpha,\beta}(x)|^4}.
\end{equation}
\end{lemma}

\begin{proof}
At a normalized base point $(0,0,x_3)$, $|x_3|<1$, we have
\[
 K=1+|x_3|^2,\qquad \Delta=(1-|x_3|^2)^2,\qquad V=1.
\]
For $\lambda_1,\lambda_2\in\mathbb D$,
\[
 |1+x_3 \lambda_1\lambda_2|\ge1-|\lambda_1||\lambda_2|
 \ge\sqrt{(1-|\lambda_1|^2)(1-|\lambda_2|^2)},
\]
where the second inequality follows by squaring and using
$(|\lambda_1|-|\lambda_2|)^2\ge0$. Equality in the supremum for $\kappa$ is attained at
$\lambda_1=\lambda_2=0$. Thus the squared fibre radius is $1=V$ at normalized points.

The normal form and \eqref{hex:K-covariance}--\eqref{hex:Delta-covariance}
prove positivity of $K$ and $\Delta$ everywhere on $\mathbb E$. In fact,
\begin{equation}\label{hex:E-explicit}
 \mathbb E=\{x:|x_3|<1,\ K(x)>2|q(x)|\}.
\end{equation}
Necessity of the second condition follows from the normal form and
covariance. Necessity of $|x_3|<1$ follows by factoring the polynomial
$h_x(\lambda,\lambda)$, which has no zeros on the closed disc.
For the converse, choose $b\in\mathbb C$ with $b^2=q(x)$ and let
\[
 A=\begin{pmatrix}x_1&b\\b&x_2\end{pmatrix}.
\]
Then $\det A=x_3$ and
$\|A\|_{\rm HS}^2=|x_1|^2+|x_2|^2+2|q|<1+|x_3|^2$.
If $\sigma_1,\sigma_2$ are the singular values of $A$, these inequalities
say
\[
 \sigma_1^2\sigma_2^2=|x_3|^2<1,
 \qquad (1-\sigma_1^2)(1-\sigma_2^2)>0.
\]
Both singular values are therefore less than $1$. Consequently
$I-A\operatorname{diag}(\lambda_1,\lambda_2)$ is invertible for $|\lambda_1|,|\lambda_2|\le1$, and its
determinant is $h_x(\lambda_1,\lambda_2)$. Thus $x\in\mathbb E$, proving
\eqref{hex:E-explicit}.

The squared actual fibre radius and $V$ both transform by the positive factor
$c/|d_{\alpha,\beta}|^2$. Their equality at normalized points therefore
proves the radius formula everywhere. Formula \eqref{hex:F-covariance}
now follows from \eqref{hex:projective-full-map} and the covariances of
$\Delta$ and $V$.

We verify explicitly the assertion about the boundary of the base. The
stable-polynomial definition of $\mathbb E$ is equivalent to
\begin{equation}\label{hex:tetrablock-scalar-criterion}
 |x_2|<1,\qquad
 B(x):=1-|x_1|^2+|x_2|^2-|x_3|^2
       >2|x_2-x_3\bar x_1|.
\end{equation}
Indeed, for $|x_2|<1$, the function
$(x_1-x_3\lambda_2)/(1-x_2\lambda_2)$ is holomorphic on a neighborhood of
$\overline{\mathbb D}$. Its modulus is strictly less than $1$ there if
and only if it is so on $|\lambda_2|=1$, by the maximum modulus principle.
The latter condition is precisely \eqref{hex:tetrablock-scalar-criterion},
by expanding $|1-x_2\lambda_2|^2-|x_1-x_3\lambda_2|^2$ on the unit circle. This modulus
condition is equivalent to the nonvanishing of $h_x(\lambda_1,\lambda_2)$ on the closed
bidisc.
A direct expansion gives
\begin{equation}\label{hex:Delta-second-form}
 \Delta(x)=B(x)^2-4|x_2-x_3\bar x_1|^2.
\end{equation}
For $x\in\overline{\mathbb E}$, we have $|x_2|\le1$ and
$B(x)\ge2|x_2-x_3\bar x_1|$. If $|x_2|<1$ and $\Delta(x)>0$,
\eqref{hex:tetrablock-scalar-criterion} implies $x\in\mathbb E$.
If $|x_2|=1$, continuity of the inequality
$|x_1-x_3\lambda_2|\le|1-x_2\lambda_2|$ at $\lambda_2=\bar x_2$ gives
$x_1=x_3\bar x_2$, and \eqref{hex:Delta-second-form} gives $\Delta(x)=0$.
It follows that
\[
 \Delta=0\quad\text{on }\partial\mathbb E.
\]

The domain $\mathbb E$ is bounded: the restrictions $h_x(\lambda_1,0)$ and
$h_x(0,\lambda_2)$ give $|x_1|,|x_2|<1$, and factoring the stable polynomial
$h_x(\lambda,\lambda)$ gives $|x_3|<1$. Also $|a|<1$ on $\mathbb H$ by
using $\lambda_1=\lambda_2=0$ in \eqref{hex:def}. Thus $\mathbb H$ is bounded.
The formulae for $\Delta$ and $V$ extend continuously to
$\overline{\mathbb E}$, with $\Delta\ge0$. At a boundary point of
$\mathbb H$ lying over an interior base point, the radius formula gives
$V=|a|^2$; over a boundary point of the base, $\Delta=0$. In both cases
\eqref{hex:F} tends to zero. This proves the asserted continuous extension
and boundary vanishing. Positivity and real analyticity in $\mathbb H$
are immediate from the same formula.
\end{proof}

\section{The strong linear convexity calculation}

For a real $C^2$ function $\rho$ define
\[
 \mathcal L_\rho(X)=\sum_{j,k}\rho_{z_j\bar z_k}X_j\bar X_k,
 \qquad
 \mathcal H_\rho(X)=\sum_{j,k}\rho_{z_jz_k}X_jX_k.
\]
Strong linear convexity at a regular boundary point is the inequality
$\mathcal L_\rho(X)>|\mathcal H_\rho(X)|$ for every nonzero complex
tangent vector $X$.

We also use that projective linear maps preserve this inequality.
Indeed, if $T(z)=(Az+b)/(c\cdot z+d)$, then
\[
 T''(z)(X,X)=-\frac{2(c\cdot X)}{c\cdot z+d}\,T'(z)X.
\]
For a complex tangent vector, the contribution of this second derivative
to the pure complex Hessian of a composed defining function vanishes.
The Levi form transforms tensorially, so the strong linear convexity
inequality is preserved.

\begin{lemma}\label{hex:hessian-lemma}
For every $y\in\mathbb H$, the function $\rho=F(y)-F$ satisfies
\[
 \mathcal L_\rho(X)>|\mathcal H_\rho(X)|
 \qquad\text{if }X\ne0\text{ and }\partial\rho(y)X=0.
\]
The assertion is also valid if $dF(y)=0$, in which case every vector is
complex tangent.
\end{lemma}

\begin{proof}
Apply Lemma~\ref{hex:automorphisms-normal-form} to send $y$ to
\[
 y_0=(A,0,0,r),\qquad 0\le A,r<1.
\]
Write
\[
 B=1-A^2>0,\qquad d=1-r^2>0,
 \qquad F(y_0)=dB.
\]
Let $T$ be this projective automorphism and let $S=T^{-1}$.
The covariance formula, also after the coordinate rotations, has the form
\[
 F(S(z))=\frac{\chi F(z)}{|E(z)|^4},
 \qquad \chi>0,
\]
where $E$ is affine holomorphic and nonzero near $y_0$. Thus multiplication
of $F(y)-F(S(z))$ by $|E(z)|^4/\chi$ gives the level function
\begin{equation}\label{hex:normalized-defining-function}
 \widetilde\rho(z)=dB\,|D(z)|^4-F(z),
 \qquad D(z)=\frac{E(z)}{E(y_0)}.
\end{equation}
We used $F(y)|E(y_0)|^4/\chi=F(y_0)=dB$.
In particular, $D$ is affine holomorphic and $D(y_0)=1$.
Multiplication of a zero-level function by a positive smooth factor
multiplies both complex Hessians on its complex kernel by that factor.
This assertion remains valid at a critical point, since all the terms
involving the first differential then vanish. Together with projective
invariance, this reduces the proof to
\eqref{hex:normalized-defining-function}, without any regularity
assumption on the level.

Use coordinates $(a,x_1,x_2,x_3)$ near $y_0$ and a vector
\[
 X=(\xi,X_1,X_2,X_3).
\]
The expansion in the middle variables is
\[
 \begin{split}
 \Delta(x_1,x_2,x_3)
 ={}&(1-|x_3|^2)^2
 -2(1+|x_3|^2)(|x_1|^2+|x_2|^2)
 +8\operatorname{Re}(\bar x_3 x_1x_2)\\
 &+(|x_1|^2-|x_2|^2)^2.
 \end{split}
\]
Consequently, the terms of order at most two in $(x_1,x_2)$ in $F$ are
\begin{equation}\label{hex:F-middle-expansion}
 \begin{split}
 F(a,x_1,x_2,x_3)={}&(1-|a|^2)(1-|x_3|^2)\\
 &-\left(1+\frac{(1-|a|^2)(1+|x_3|^2)}{1-|x_3|^2}\right)
       (|x_1|^2+|x_2|^2)\\
 &+2\left(1+\frac{2(1-|a|^2)}{1-|x_3|^2}\right)
       \operatorname{Re}(\bar x_3 x_1x_2)
 +O((|x_1|+|x_2|)^4).
 \end{split}
\end{equation}
This expansion is valid locally with $|x_3|<1$.
At $y_0$ we therefore have
\[
 \partial F(y_0)X=-dA\xi-BrX_3.
\]
Writing $\ell(X)=D'(y_0)X$, complex tangency to
\eqref{hex:normalized-defining-function} is
\[
 2dB\,\ell(X)+dA\xi+BrX_3=0.
\]
Thus
\begin{equation}\label{hex:tangent-denominator}
 \ell(X)=-\frac{dA\xi+BrX_3}{2dB}.
\end{equation}
Since at $D=1$ the Levi and pure complex Hessians of $|D|^4$ are
$4|\ell(X)|^2$ and $2\ell(X)^2$, respectively, substitution of
\eqref{hex:tangent-denominator} into the Hessians obtained from
\eqref{hex:F-middle-expansion} gives
\begin{align}
 \mathcal L_{\widetilde\rho}(X)
 &=\frac dB|\xi|^2+\frac Bd|X_3|^2
   +\left(1+\frac{B(1+r^2)}d\right)(|X_1|^2+|X_2|^2),
 \label{hex:Levi-normal-form}\\
 \mathcal H_{\widetilde\rho}(X)
 &=\frac12\left(\sqrt{\frac dB}A\xi
                -\sqrt{\frac Bd}rX_3\right)^2
   -2r\left(1+\frac{2B}d\right)X_1X_2.
 \label{hex:pure-normal-form}
\end{align}
For clarity, the cancellation in the first formula comes from the
opposite mixed terms $2Ar\operatorname{Re}(\xi\bar X_3)$ contributed by
$|D|^4$ and $F$.

Put $L_0=\frac dB|\xi|^2+\frac Bd|X_3|^2$. Cauchy--Schwarz and
$2|X_1X_2|\le|X_1|^2+|X_2|^2$ yield
\[
 |\mathcal H_{\widetilde\rho}(X)|
 \le\frac{A^2+r^2}{2}L_0
    +r\left(1+\frac{2B}d\right)(|X_1|^2+|X_2|^2).
\]
Combining this with \eqref{hex:Levi-normal-form},
\begin{equation}\label{hex:strict-positive-bound}
 \begin{aligned}
 \mathcal L_{\widetilde\rho}(X)
 -|\mathcal H_{\widetilde\rho}(X)|
 \ge{}&\left(1-\frac{A^2+r^2}{2}\right)
       \left(\frac dB|\xi|^2+\frac Bd|X_3|^2\right)\\
 &+\left(1-r+\frac{B(1-r)^2}{d}\right)
       (|X_1|^2+|X_2|^2)>0
 \end{aligned}
\end{equation}
for $X\ne0$. Both coefficients are strictly positive because $A,r<1$.
This proves the lemma.
\end{proof}

\section{Regularity, connectedness, and exhaustion}
\label{sec:regularity-exhaustion}

\begin{proof}[Proof of Theorem~\ref{hex:exhaustion-theorem}]
The radius assertion and the boundary properties of $F$ have already been
proved. We first show that the origin is the only critical point of $F$ in
$\mathbb H$. The origin is indeed critical, by the formula for $F$, and
$F(0)=1$.

The function $F$ is invariant under the weighted circle action
\[
 (a,x_1,x_2,x_3)\longmapsto
 (e^{it}a,e^{it}x_1,e^{it}x_2,e^{2it}x_3).
\]
Suppose that $dF(y)=0$ at a nonzero point $y\in\mathbb H$. Apply
Lemma~\ref{hex:hessian-lemma} with $\rho=F(y)-F$. Every vector is complex
tangent at $y$, so the inequality in that lemma implies that the real
Hessian of $\rho$ is positive definite at $y$. On the other hand, the
weighted circle orbit through $y$ is nonconstant and $\rho$ is constant
along it. Its tangent vector is
\[
 Z=i(a,x_1,x_2,2x_3)\ne0.
\]
Differentiating twice along the orbit, and using $d\rho(y)=0$, gives
$\operatorname{Hess}_{\mathbb R}\rho(y)[Z,Z]=0$, a contradiction.
Thus $0$ is the unique critical point.

Since $\mathbb H$ is bounded and $F$ extends continuously by zero to its
boundary, $F$ attains its maximum on $\overline{\mathbb H}$. Its
positivity inside forces a maximum point to be interior. The uniqueness
of the critical point therefore gives
\[
 0<F\le1\quad\text{on }\mathbb H,
 \qquad F^{-1}(1)=\{0\}.
\]
For $0<\varepsilon<1$, the set $\{F\ge\varepsilon\}$ is consequently
a compact subset of $\mathbb H$, and $dF\ne0$ on $\{F=\varepsilon\}$.
Thus $D_\varepsilon$ has real-analytic regular boundary. The Hessian
lemma proves its strong linear convexity at every boundary point.

It remains to verify connectedness rather than merely selecting a
component. Let $\mathcal C$ be a connected component of $\{F>\varepsilon\}$.
Its closure is compact in $\mathbb H$, and $F=\varepsilon$ on
$\partial\mathcal C$. Since $\mathcal C$ contains points where $F>\varepsilon$, the
maximum of $F$ on $\overline{\mathcal C}$ occurs in $\mathcal C$. This is a critical point,
so $0\in\mathcal C$. Every component therefore contains the origin, and there
can be only one component.

Finally, if $0<\varepsilon'<\varepsilon<1$, then
\[
 \overline{D_\varepsilon}
 =\{F\ge\varepsilon\}
 \subset\{F>\varepsilon'\}=D_{\varepsilon'},
\]
with compact containment. Every point of $\mathbb H$ has $F>0$, so it
belongs to $D_\varepsilon$ for all sufficiently small positive
$\varepsilon$. This proves the asserted exhaustion.
\end{proof}

\section{Hyperconvexity and tautness}
\label{sec:hyperconvexity}

We shall need tautness only to formulate the conclusion in the usual
Lempert-domain form.  It follows directly from the defining family of the
hexablock.

For $\lambda_2\in\overline{\mathbb D}$ and
$\lambda=(\lambda_1,\lambda_2)\in\overline{\mathbb D}^{\,2}$ put
\begin{align}
 \phi_{\lambda_2}(x)
 &=\frac{x_1-x_3\lambda_2}{1-x_2\lambda_2},\label{eq:phi-test}\\
 \psi_{\lambda}(a,x)
 &=\frac{a\sqrt{(1-|\lambda_1|^2)(1-|\lambda_2|^2)}}
 {1-x_1\lambda_1-x_2\lambda_2+x_3\lambda_1\lambda_2}.
 \label{eq:psi-test}
\end{align}
The denominator in \eqref{eq:phi-test} does not vanish on $\mathbb E$,
and
\[
 h_x(\lambda_1,\lambda_2)
 =(1-x_2\lambda_2)\bigl(1-\lambda_1\phi_{\lambda_2}(x)\bigr).
\]
Hence $\phi_{\lambda_2}\in\mathcal O(\mathbb E,\mathbb D)$.  Likewise,
by the definition of $\mathbb H$, $\psi_\lambda\in\mathcal O(\mathbb H,\mathbb D)$.

Define
\begin{align*}
 \mathcal A(x)
 &=\max_{\lambda_2\in\overline{\mathbb D}}
       |\phi_{\lambda_2}(x)|,\\
 \mathcal B(a,x)
 &=\max_{\lambda\in\overline{\mathbb D}^{\,2}}
       |\psi_\lambda(a,x)|,
\end{align*}
and
\begin{equation}\label{eq:U-psh}
 U(a,x)=\max\{|x_2|,\mathcal A(x),\mathcal B(a,x)\}-1.
\end{equation}

\begin{proposition}\label{prop:hyperconvex}
The function $U$ is a continuous negative plurisubharmonic exhaustion of
$\mathbb H$.  In particular, $\mathbb H$ is hyperconvex and taut.
\end{proposition}

\begin{proof}
For fixed parameters, $|\phi_{\lambda_2}|$ and $|\psi_\lambda|$ are
plurisubharmonic.  On every relatively compact subset of the corresponding
domain the denominators are uniformly separated from zero, uniformly in the
compact parameter sets.  It follows that $\mathcal A$ and $\mathcal B$ are
continuous.  Since they are continuous suprema of locally uniformly bounded
families of plurisubharmonic functions, they are plurisubharmonic.  Thus $U$
is continuous and plurisubharmonic.  Moreover,
\[
 -1\le U<0\qquad\text{on }\mathbb H.
\]

It remains to prove that $U$ tends to $0$ at the boundary.  Let
$(a^j,x^j)\to(a^0,x^0)\in\partial\mathbb H$ with $(a^j,x^j)\in\mathbb H$.
If $x^0\in\mathbb E$, then the fibre inequality must become sharp, and hence
$\mathcal B(a^j,x^j)\to1$.

Suppose that $x^0\in\partial\mathbb E$.  If $|x_2^0|=1$, then
$|x_2^j|\to1$.  If $|x_2^0|<1$, then, because $x^0$ is not an interior point
of the stable-polynomial domain $\mathbb E$, there is
$\lambda^0=(\lambda_1^0,\lambda_2^0)\in\overline{\mathbb D}^{\,2}$ such that
$h_{x^0}(\lambda_1^0,\lambda_2^0)=0$.  Since
$1-x_2^0\lambda_2^0\ne0$, the factorization above gives
\[
 \lambda_1^0\phi_{\lambda_2^0}(x^0)=1.
\]
The function $\phi_{\lambda_2^0}$ is holomorphic near $x^0$, while
$|\phi_{\lambda_2^0}(x^j)|<1$.  Hence
$|\phi_{\lambda_2^0}(x^0)|\le1$, and the displayed identity forces equality.
Thus $\mathcal A(x^j)\to1$.  In all cases,
\[
 U(a^j,x^j)\longrightarrow0.
\]
Since $\mathbb H$ is bounded, the sublevel sets $\{U\le c\}$, $c<0$, are
relatively compact in $\mathbb H$.  This proves hyperconvexity.

For completeness, we give the standard normal-family argument for tautness.
Let $f_j\in\mathcal O(\mathbb D,\mathbb H)$ be a sequence which is not
compactly divergent.  Then there exist compact sets $K\Subset\mathbb D$ and
$L\Subset\mathbb H$ and, after passing to a subsequence, points
$\lambda_j\in K$ such that $f_j(\lambda_j)\in L$.  Passing to a further
subsequence, we may suppose that
\[
 \lambda_j\longrightarrow\lambda_0\in K,
 \qquad
 f_j(\lambda_j)\longrightarrow y_0\in L.
\]
Since $\mathbb H$ is bounded, Montel's theorem gives, after another passage
to a subsequence, local uniform convergence to a holomorphic map
$f:\mathbb D\to\overline{\mathbb H}$.  The local uniform convergence and
$\lambda_j\to\lambda_0$ give $f(\lambda_0)=y_0\in\mathbb H$.  Put
$u_j=U\circ f_j$.  The functions $u_j$ are subharmonic and satisfy
$-1\le u_j<0$.  Hence
\[
 u=(\limsup_{j\to\infty}u_j)^*
\]
is subharmonic.  Since $f(\lambda_0)$ is an interior point, $f$ maps a
neighborhood of $\lambda_0$ into a compact subset of $\mathbb H$; there
$u_j\to U\circ f$ locally uniformly.  Thus $u(\lambda_0)<0$.
If $f(\lambda_*)\in\partial\mathbb H$ for some $\lambda_*\in\mathbb D$,
then the exhaustion property of $U$ gives $u_j(\lambda_*)\to0$.  Since
$u\le0$, we obtain $u(\lambda_*)=0$, contradicting the maximum principle.
Therefore $f(\mathbb D)\subset\mathbb H$, proving tautness.
\end{proof}

\section{The tetrablock, pentablock and symmetrized bidisc as sections}
\label{sec:sections}

The purpose of this section is to extract all lower-dimensional consequences
without repeating the strong-linear-convexity calculation.

\subsection{The three identifications}

Let
\[
 S_{\mathbb E}=\{(a,x):a=0\},\qquad
 S_{\mathbb P}=\{(a,x):x_1=x_2\},
\]
\[
 S_{\mathbb G}=S_{\mathbb E}\cap S_{\mathbb P}
               =\{(a,x):a=0,\ x_1=x_2\}.
\]
The tetrablock identification is immediate:
\begin{equation}\label{eq:E-section}
 \iota_{\mathbb E}:\mathbb E\longrightarrow \mathbb H\cap S_{\mathbb E},
 \qquad
 \iota_{\mathbb E}(x_1,x_2,x_3)=(0,x_1,x_2,x_3),
\end{equation}
is a complex linear biholomorphism.

For the pentablock, define
\begin{equation}\label{eq:P-section-map}
 \iota_{\mathbb P}(a,s,p)=\left(a,\frac{s}{2},\frac{s}{2},p\right).
\end{equation}
We verify that this identifies $\mathbb P$ with $\mathbb H\cap S_{\mathbb P}$.
We first record the corresponding identification of the bases:
\begin{equation}\label{eq:G-E-diagonal}
 (s,p)\in\mathbb G_2
 \quad\Longleftrightarrow\quad
 \left(\frac{s}{2},\frac{s}{2},p\right)\in\mathbb E.
\end{equation}
Indeed, suppose first that
$x=(s/2,s/2,p)\in\mathbb E$.  Then
\[
 h_x(t,t)=1-st+pt^2
\]
has no zero in $\overline{\mathbb D}$.  If $\lambda_1,\lambda_2$ are the
roots, counted with multiplicity, of
\[
 \zeta^2-s\zeta+p=0,
\]
then
\[
 1-st+pt^2=(1-\lambda_1t)(1-\lambda_2t).
\]
If $|\lambda_j|\ge1$ for some $j$, then $t=1/\lambda_j$ is a zero in
$\overline{\mathbb D}$, a contradiction.  Thus
$\lambda_1,\lambda_2\in\mathbb D$ and
$(s,p)=(\lambda_1+\lambda_2,\lambda_1\lambda_2)\in\mathbb G_2$.

Conversely, let $(s,p)=\pi(\lambda)$ with
$\lambda=(\lambda_1,\lambda_2)\in\mathbb D^2$, and put
\[
 U_\lambda=|1-\overline{\lambda_2}\lambda_1|,
 \qquad
 W_\lambda=\sqrt{(1-|\lambda_1|^2)(1-|\lambda_2|^2)}.
\]
For $x=(s/2,s/2,p)$ one has
\begin{align*}
 K(x)&=1-\frac{|s|^2}{2}+|p|^2
      =\frac{U_\lambda^2+W_\lambda^2}{2},\\
 4|q(x)|&=|\lambda_1-\lambda_2|^2
          =U_\lambda^2-W_\lambda^2.
\end{align*}
Hence
\[
 K(x)-2|q(x)|=W_\lambda^2>0,
 \qquad
 |x_3|=|p|<1.
\]
The characterization \eqref{hex:E-explicit} therefore gives $x\in\mathbb E$, and
\eqref{eq:G-E-diagonal} follows.  The same identities also give
\[
 \sqrt{\Delta(x)}=U_\lambda W_\lambda,
\]
and consequently
\begin{equation}\label{eq:V-penta}
 V\left(\frac{s}{2},\frac{s}{2},p\right)
 =\frac{(U_\lambda+W_\lambda)^2}{4}.
\end{equation}
The fibre-radius formula of Agler--Lykova--Young
\cite[Theorem~5.2]{ALY2015} says precisely that, for
$(s,p)=\pi(\lambda)\in\mathbb G_2$,
\[
 (a,s,p)\in\mathbb P
 \quad\Longleftrightarrow\quad
 |a|<\frac{U_\lambda+W_\lambda}{2}.
\]
Together with \eqref{eq:G-E-diagonal}, \eqref{eq:V-penta}, and the radius
formula for $\mathbb H$, this proves
\begin{equation}\label{eq:P-section}
 \iota_{\mathbb P}(\mathbb P)=\mathbb H\cap S_{\mathbb P}.
\end{equation}

Finally, the zero fibre of the pentablock is $\{0\}\times\mathbb G_2$.
Consequently
\begin{equation}\label{eq:G-section}
 \iota_{\mathbb G}(s,p)=\left(0,\frac{s}{2},\frac{s}{2},p\right)
\end{equation}
identifies $\mathbb G_2$ with $\mathbb H\cap S_{\mathbb G}$.

\subsection{Restriction of the single exhaustion}

\begin{proposition}\label{prop:section-exhaustions}
For $S\in\{S_{\mathbb E},S_{\mathbb P},S_{\mathbb G}\}$ and
$0<\varepsilon<1$, the set
\[
 D_\varepsilon^S=D_\varepsilon\cap S
\]
is a connected domain in $S$ with real-analytic strongly linearly convex
boundary.  The family $D_\varepsilon^S$ is compactly nested and exhausts
$\mathbb H\cap S$ as $\varepsilon\downarrow0$.
\end{proposition}

\begin{proof}
There are two elementary symmetries of $F$ which are exactly adapted to the
sections.  First, $F$ depends on $a$ only through $|a|^2$, and hence at
$a=0$ its differential vanishes in the complex normal $a$-direction.
Second, $F$ is invariant under the transposition
\[
 (a,x_1,x_2,x_3)\longmapsto(a,x_2,x_1,x_3).
\]
Hence, at a point with $x_1=x_2$, its differential vanishes on the complex
normal line
\[
 \{(0,\zeta,-\zeta,0):\zeta\in\mathbb C\}.
\]
Therefore a critical point of $F$ restricted to any of the three sections
is a critical point of $F$ on $\mathbb H$.  By the proof of
Theorem~\ref{hex:exhaustion-theorem}, the only such point is the origin,
where $F=1$.  Thus every restricted level $F=\varepsilon$, $0<\varepsilon<1$,
is regular.

If $X\ne0$ is complex tangent to a restricted level, then it is also complex
tangent to the full level.  Lemma~\ref{hex:hessian-lemma} therefore gives
\[
 \mathcal L_{\varepsilon-F}(X)>
 |\mathcal H_{\varepsilon-F}(X)|,
\]
so the restricted level is strongly linearly convex.  Connectedness follows
from exactly the same maximum argument as for $D_\varepsilon$: every
component of $\{F|_S>\varepsilon\}$ contains an interior maximum of $F|_S$,
hence contains the origin.  Compact nesting and exhaustion are inherited
from the corresponding properties of $D_\varepsilon$.
\end{proof}

Transporting these domains through
\eqref{eq:E-section}, \eqref{eq:P-section-map}, and \eqref{eq:G-section}
gives explicit exhaustions
\begin{align*}
 E_\varepsilon
 &=\{x\in\mathbb E:F(0,x)>\varepsilon\},\\
 P_\varepsilon
 &=\left\{(a,s,p)\in\mathbb P:
     F\left(a,\frac{s}{2},\frac{s}{2},p\right)>\varepsilon\right\},\\
 G_\varepsilon
 &=\left\{(s,p)\in\mathbb G_2:
     F\left(0,\frac{s}{2},\frac{s}{2},p\right)>\varepsilon\right\}.
\end{align*}
There is no additional Hessian calculation for any of these three families.

\begin{corollary}\label{cor:sections-hyperconvex}
The domains $\mathbb E$, $\mathbb P$ and $\mathbb G_2$ are hyperconvex and
taut.
\end{corollary}

\begin{proof}
Restrict the continuous negative plurisubharmonic exhaustion $U$ from
Proposition~\ref{prop:hyperconvex} to the corresponding closed complex
linear section and transport it by the linear identifications above.  The
restriction is still negative and plurisubharmonic, and its sublevel sets are
relatively compact in the section.
\end{proof}

\section{Passage to invariant distances and metrics}
\label{sec:metrics}

We isolate the standard exhaustion argument in a form that applies
simultaneously to all four domains.

\begin{proposition}\label{prop:metric-limit}
Let $n\ge2$ and let $\Omega\subset\mathbb C^n$ be a domain admitting an increasing
exhaustion
\[
 \Omega_1\Subset\Omega_2\Subset\cdots,\qquad
 \bigcup_{j=1}^\infty\Omega_j=\Omega,
\]
where every $\Omega_j$ is bounded, has $C^2$ boundary, and is strongly
linearly convex.  Then
\[
 c_\Omega=k_\Omega=\ell_\Omega,
 \qquad
 \gamma_\Omega=\kappa_\Omega.
\]
\end{proposition}

\begin{proof}
By the Lempert theorem in the form of Kosi\'nski--Warszawski
\cite[Theorem~1.3]{KW2013},
\[
 c_{\Omega_j}=\ell_{\Omega_j},
 \qquad
 \gamma_{\Omega_j}=\kappa_{\Omega_j}.
\]
Fix $z,w\in\Omega$ and omit finitely many indices.  Monotonicity gives a
decreasing sequence $c_{\Omega_j}(z,w)$ bounded below by $c_\Omega(z,w)$.
Choose approximate Carath\'eodory extremals
$f_j:\Omega_j\to\mathbb D$, normalized by $f_j(z)=0$.  A diagonal Montel
argument gives a subsequence converging locally uniformly on $\Omega$ to a
holomorphic function $f:\Omega\to\overline{\mathbb D}$.  Since $f(z)=0$,
the maximum principle gives $f(\Omega)\subset\mathbb D$.  Evaluation at $w$
shows
\[
 c_{\Omega_j}(z,w)\downarrow c_\Omega(z,w).
\]
Consequently
\[
 \ell_\Omega(z,w)
 \le \ell_{\Omega_j}(z,w)
 =c_{\Omega_j}(z,w)\longrightarrow c_\Omega(z,w).
\]
Together with $c_\Omega\le k_\Omega\le\ell_\Omega$, this proves the distance
equalities.

The infinitesimal argument is identical.  Normalized approximate extremals
for $\gamma_{\Omega_j}(z;X)$ have a locally uniformly convergent subsequence,
and Cauchy's formula gives convergence of the derivatives.  Hence
\[
 \gamma_{\Omega_j}(z;X)\downarrow\gamma_\Omega(z;X).
\]
By monotonicity,
\[
 \kappa_\Omega(z;X)
 \le\kappa_{\Omega_j}(z;X)
 =\gamma_{\Omega_j}(z;X)
 \longrightarrow\gamma_\Omega(z;X).
\]
The reverse inequality $\gamma_\Omega\le\kappa_\Omega$ completes the proof.
\end{proof}

\begin{proof}[Proof of Corollary~\ref{cor:intro-domains}]
Apply Proposition~\ref{prop:metric-limit} first to the exhaustion
$D_{2^{-j}}$ of $\mathbb H$, and then to its three restrictions from
Proposition~\ref{prop:section-exhaustions}.  Tautness of $\mathbb H$ is
Proposition~\ref{prop:hyperconvex}; tautness of the three sections is
Corollary~\ref{cor:sections-hyperconvex}.  Hence all four domains are
Lempert domains in the sense fixed in the introduction.

Strong linear convexity implies $\mathbb C$-convexity.  Since the
exhaustions are increasing, the same is true for their unions.  Indeed, fix a
complex affine line $L$ and a loop in $L\cap\Omega$.  Its image is compact,
hence it is contained in $L\cap\Omega_j$ for some sufficiently large $j$.
This intersection is nonempty and simply connected, so the loop contracts
there.  Thus $\mathbb H$ is $\mathbb C$-convex, and the same conclusion holds
for the three sections.
\end{proof}

\medskip
\noindent\textbf{Use of artificial intelligence.}
ChatGPT (OpenAI) was used in the preparation of this manuscript for editorial
assistance and for carrying out calculations related to the construction of the exhausting families considered in the paper.

\end{document}